\documentclass[12pt]{article}
\usepackage{amsmath,amsfonts,cite}

\newtheorem{thm}{Theorem}[section]
\newtheorem{lem}[thm]{Lemma}
\newtheorem{cor}[thm]{Corollary}

\def\pf{\noindent{\it Proof.} }
\def\qed{\nopagebreak\hfill{\rule{4pt}{7pt}}\medbreak}

\makeatletter \@addtoreset{equation}{section} \makeatother

\begin{document}

\begin{center}
{\Large\bf On Cui-Kano's Characterization Problem

on Graph Factors}
\end{center}

\begin{center}
H. L. Lu$^{1}$ and David G. L. Wang$^{2}$\\[6pt]

$^{1}$Department of Mathematics\\
Xi'an Jiaotong University, Xi'an 710049, P. R. China\\
{\tt $^{1}$luhongliang215@sina.com}

$^{2}$Beijing International Center for Mathematical Research\\
Peking University, Beijing 100871, P. R. China\\
{\tt $^{2}$wgl@math.pku.edu.cn}
\end{center}

\begin{abstract}
An $H_n$-factor of a graph~$G$
is defined to be a spanning subgraph~$F$ of~$G$
such that each vertex has degree
belonging to the set $\{1,3,5,\ldots,2n-1,2n\}$ in~$F$.
In this paper,
we investigate $H_n$-factors of graphs
by using Lov\'asz's structural descriptions to
the degree prescribed subgraph problem.
We find some sufficient conditions for
the existence of an $H_n$-factor of a graph.
In particular,
we make progress on the characterization problem for
a special family of graphs
proposed by Cui and Kano in~1988.
\end{abstract}


\noindent\textbf{2010 AMS Classification:} 05C75

\section{Introduction}

Let~$G$~be a simple graph with vertex set~$V(G)$.
Let $f,g\colon V(G)\to\mathbb{Z}$ be functions on~$V(G)$.
A $(g,f)$-factor of~$G$ is a spanning subgraph $F$ such that
\[
g(v)\le d_F(v)\le f(v)
\]
for any vertex $v$,
where $d_F(v)$ is the degree of~$v$ in~$F$.
In particular, if there exist integers~$a$ and~$b$ such that
$g(v)=a$ and $f(v)=b$ for all vertices~$v$,
then the $(g,f)$-factor is called an $[a,b]$-factor.
For example,
connected $[2,2]$-factors are nothing but Hamiltonian cycles,
and $[1,1]$-factors are perfect matchings.
There is a large amount of literature on graph factors,
see Plummer~\cite{Plu07},
Liu and Yu~\cite{YL09},
and Akiyama and Kano~\cite{AK11} for surveys.
For connected factors,
we refer the reader to Kouider and Vestergaard~\cite{KV05}.

Let~$H$~be a function associating
a subset of $\mathbb{Z}$ with each vertex of~$G$,
called a degree prescription.
It is natural to generalize $(g,f)$-factors to $H$-factors,
i.e., spanning subgraphs~$F$ such that
\begin{equation}\label{cond_fac}
d_F(v)\in H(v)
\end{equation}
for all vertices~$v$. Let~$F$~be a spanning subgraph of~$G$.
Following Lov\'asz~\cite{Lov72}, one may measure the ``deviation''
of~$F$ from the condition~\eqref{cond_fac} by
\begin{equation}\label{eq_deviation}
\delta_H(F)=\sum_{v\in V(G)}
\min\bigl\{|d_F(v)-h|\,\colon\,h\in H(v)\bigr\}.
\end{equation}
Moreover, the ``solvability'' of~\eqref{cond_fac}
can be characterized by
\[
\delta(H)=\min\{\delta_H(F)\,\colon\,\mbox{$F$ is a spanning subgraph of~$G$}\}.
\]
The subgraph $F$ is said to be $H$-optimal if
$\delta_H(F)=\delta(H)$.
It is clear that
$F$ is an $H$-factor if and only if $\delta_H(F)=0$,
and any $H$-factor (if exists) is $H$-optimal.

In~\cite{Lov69},
Lov\'asz proposed the problem of
determining the value of $\delta(H)$,
called the degree prescribed subgraph problem.
Let
\[
H=\{h_1,h_2,\ldots,h_n\}
\]
be a set of integers, where $h_1<h_2<\cdots<h_n$.
It is said to be an allowed set if each of its gaps
has at most one integer, i.e.,
\[
h_{i+1}-h_i\le2,\quad\forall\, 1\le i\le n-1.
\]
We say that a prescription $H$ is allowed if
$H(v)$ is an allowed set for all vertices $v$.
Lov\'asz~\cite{Lov72} built up a whole theory
to the degree prescribed subgraph problem
in case that $H$ is an allowed prescription.
He showed that the problem is NP-complete
without the restriction that $H$ is allowed.
Cornu\'ejols~\cite{Cor88} provided the first polynomial algorithm
for the problem with~$H$ allowed.

A special case of the degree prescribed subgraph problem
is the so-called $f$-parity subgraph problem, i.e.,
the problem with
\[
H(v)=\{\ldots,\,f(v)-4,\,f(v)-2,\,f(v)\}
\]
for some function $f\colon V(G)\to\mathbb{Z}$.
The first investigation of the $f$-parity subgraph problem
is due to Amahashi~\cite{Ama85},
who gave a Tutte type characterization for graphs
having a global odd factor.
Let~$S$~be a subset of~$V(G)$.
Denote by~$G-S$ the subgraph of~$G$
obtained by removing all vertices in $S$.
Denote by~$o(G)$ the number of odd components of~$G$.
Let $n\ge2$ be an integer
independent of the number of vertices of~$G$.
Let $H_o$ be the prescription associating the first
$n$ positive odd integers with each vertex, i.e.,
\[
H_o(v)=\{1,3,5,\ldots,2n-1\}.
\]
\begin{thm}[Amahashi]\label{thm_Amahashi}
A graph $G$ has an $H_o$-factor if and only if
\begin{equation}\label{cond_Amahashi}
o(G-S)\le(2n-1)\,|S|,\quad\forall\, S\subseteq V(G).
\end{equation}
\end{thm}

For general odd value functions~$f$,
Cui and Kano~\cite{CK88} established a Tutte type theorem.
Noticing the form of the condition~\eqref{cond_Amahashi},
they asked the question of characterizing graphs~$G$
in terms of graph factors such that
\begin{equation}\label{cond_CK}
o(G-S)\le 2n\,|S|,\quad\forall\, S\subseteq V(G).
\end{equation}
For more studies on the $f$-parity subgraph problem,
see Topp and Vestergaard~\cite{TV93},
and Kano, Katona, and Szab\'o~\cite{KKS09}.

Motivating by solving Cui-Kano's problem,
we consider the degree prescribed subgraph problem
for the special prescription
\begin{equation}\label{Hn}
H_n(v)=H_o(v)\cup\{2n\}=\{1,3,5,\ldots,2n-1,2n\}.
\end{equation}
We shall study the structure of graphs which
have no $H_n$-factors
by using Lov\'asz's theory~\cite{Lov72}.
Consequently,
we obtain that any graph satisfying the condition~\eqref{cond_CK}
contains an $H_n$-factor.

Besides many applications of
Lov\'asz's structural description
for special families of graphs and special allowed prescriptions,
there is much attention paid to finding sufficient conditions
for the existence of an $H$-factor in a graph
for special prescriptions~$H$,
see~\cite{Plu07}.
In this paper,
we also give some sufficient conditions for the
existence of an $H_n$-factor in a graph.

\section{The main result}

In this section, we study $H_n$-factors of graphs
based on Lov\'asz's structural description
of the degree prescribed subgraph problem.

Let $H$ be an allowed prescription.
Denote by $I_H(v)$ the set of vertex degrees
in all $H$-optimal subgraphs, i.e.,
\[
I_H (v)=\{d_F(v)\ |\ \mbox{$F$ is $H$-optimal}\}.
\]
Comparing the set $I_H(v)$ with $H(v)$,
one may partition the vertex set $V(G)$
into four classes:
\begin{align*}
C_H&=\{v\in V(G)\,\colon\,I_H(v)\subseteq H(v)\},\\[5pt]
A_H&=\{v\in V(G)\backslash C_H\,\colon\,\min I_H(v)\ge \max H(v)\},\\[5pt]
B_H&=\{v\in V(G)\backslash C_H\,\colon\,\max I_H(v)\le \min H(v)\},\\[5pt]
D_H&=V(G)\backslash A_H\backslash B_H\backslash C_H.
\end{align*}
It is clear that
the $4$-tuple $(A_H,B_H,C_H,D_H)$
is a pairwise disjoint partition of $V(G)$.
We call it the $H$-decomposition of~$G$.
In fact, the four subsets can be distinguished
according to the contributions of their members
to the deviation~\eqref{eq_deviation}.
A graph $G$ is said to be $H$-critical if it is connected and $D_H=V(G)$.

In~\cite[Theorem~(2.1)]{Lov72},
Lov\'asz gave the following property for the subset $D_H$.

\begin{lem}[Lov\'asz]\label{lem_interval}
If $D_H\ne\emptyset$, then
the intersection
\[
[\min I_H(v),\,\max I_H(v)]\cap H(v)
\]
contains no consecutive integers for any vertex $v\in D_H$.
\end{lem}

In~\cite[Corollary (2.4)]{Lov72},
Lov\'asz gave the following structural result.

\begin{lem}[Lov\'asz]\label{lem_CD}
There is no edge between $C_H$ and $D_H$.
\end{lem}

For any subset $S\subseteq V(G)$,
denote by $G[S]$ the subgraph induced by~$S$.
Denote the number of components of $G[S]$ by $c(S)$,
and the number of odd components of $G[S]$ by $o(S)$.
In~\cite[Theorem (4.3)]{Lov72},
Lov\'asz established the formula
\[
\delta(H)=
c(D_H)+\sum_{v\in B_H}\min H(v)-
\sum_{v\in A_H}\max H(v)-\sum_{v\in B_H }d_{G-A_H}(v).
\]
By definition, $G$ contains no $H$-factors
if and only if $\delta(H)>0$.
This yields the next lemma immediately.

\begin{lem}[Lov\'asz]\label{lem_Lovasz1}
A graph $G$ contains no $H$-factors if and only if
\begin{equation}\label{ineq_lovasz1}
c(D_H)+\sum_{v\in B_H}\min H(v)>
\sum_{v\in A_H}\max H(v)+\sum_{v\in B_H }d_{G-A_H}(v).
\end{equation}
\end{lem}


Let $X\subseteq V(G)$.
For any vertex $v$, define
\begin{equation}\label{H_XY}
H_{X}(v)=\{h-e(v,X)\,\colon\,h\in H(v)\},
\end{equation}
where $e(v,X)$ is the number of edges from $v$ to $X$.
In~\cite[Theorem~(4.2)]{Lov72},
Lov\'asz showed that each component
of $G[D_H]$ is $H'$-critical where $H'=H_{B_H}$.
He~\cite[Lemma~(4.1)]{Lov72} also obtained that
$\delta(H)=1$ if $G$ is $H$-critical.
This leads to the next lemma.

\begin{lem}[Lov\'asz]\label{lem_delta=1}
If $D_H\ne\emptyset$, then
for any component~$T$ of the subgraph~$G[D_H]$,
and any $H$-optimal subgraph~$F$ of~$T$, we have
$\delta_{H'}(F)=1$.
\end{lem}

This paper concerns $H_n$-decompositions
where $H_n$ is defined by~\eqref{Hn}.
For convenience, we often use another prescription $H_n^*$
defined by
\begin{equation}\label{Hn*}
H_n^*(v)=H_n(v)\cup\{-1\}=\{-1,1,3,5,\ldots,2n-1,2n\}.
\end{equation}
Here is the main result of this paper.

\begin{thm}\label{thm_main}
Let $G$ be a graph without odd components.
If $G$ contains no $H_n$-factors,
then there exists a nonempty subset $S\subset V(G)$ such that
the subgraph $G-S$ contains at least $2n|S|+1$ odd components,
each of which contains no $H_n$-factors.
\end{thm}

\pf By the definition~\eqref{Hn*},
the graph $G$ contains no $H_n^*$-factors.
Let $(A,B,C,D)$
be the $H_n^*$-decomposition of~$G$.
We shall show that the subset~$A$
can be taken as the required~$S$.

Since $\min H_n^*=-1$,
we have $B=\emptyset$, and the inequality~\eqref{ineq_lovasz1} reduces to
\begin{equation}\label{eq:1}
c(D)>2n|A|.
\end{equation}
This implies $D\ne\emptyset$.
Let~$T$~be a component of the subgraph $G[D]$,
and~$F$~an $H_n^*$-optimal subgraph of~$T$.
Since $B=\emptyset$,
we see that $H_{B}=H_n^*$.
So $\delta_{H_n^*}(F)=1$ by Lemma~\ref{lem_delta=1}.
Therefore,
there exists a vertex, say $v_0\in V(T)$, such that
\begin{equation}\label{eq1}
\min\bigl\{\,|d_F(v)-h|\,\colon\,h\in H_n^*(v)\bigr\}
=\left\{\begin{array}{ll}
1, & \hbox{if $v=v_0$;} \\[5pt]
0, & \hbox{if $v\in V(T)\backslash\{v_0\}$.}
\end{array}
\right.
\end{equation}
On the other hand, assume that
$\max I_{H_n^*}(v)\ge 2n$ for some $v\in D$.
By Lemma~\ref{lem_interval}, we have
\[
\min I_{H_n^*}(v)\ge 2n.
\]
It follows immediately that $v\in A$,
a contradiction.
Thus $\max I_{H_n^*}(v)\le 2n-1$, namely
\[
d_F(v)\le 2n-1
\]
for any vertex $v\in T$.
Consequently, the formula~\eqref{eq1} implies that the degree
$d_{F}(v_0)$ is even, while
$d_{F}(v)$ is odd for any $v\in V(T)\backslash\{v_0\}$.
Since the sum $\sum_{v\in V(T)}d_{F}(v)$ is even,
we deduce that $T$ is an odd component of $G[D]$.

Assume that $A=\emptyset$.
Since $B=\emptyset$,
by Lemma~\ref{lem_CD},
we see that $T$ is an odd component of $G$.
But $G$ has no odd components, a contradiction.
So $A\ne\emptyset$.

By Lemma~\ref{lem_CD} and the inequality~\eqref{eq:1}, we have
\begin{equation}\label{ineq5}
2n|A|<c(D)
=o(D)
\le o(C)+o(D)
=o(C\cup D)
=o(G-A).
\end{equation}
Namely, the subgraph $G-A$
has at least $2n|A|+1$
odd components.
In view of~\eqref{eq1},
any component~$T$ of $G[D]$ has no $H_n^*$-factors.
Hence~$T$ has no $H_n$-factors.
This completes the proof. \qed

For graphs satisfying Cui-Kano's condition~\eqref{cond_CK},
we obtain the following corollary immediately.

\begin{cor}\label{cor_CK}
Any graph $G$ satisfying the condition~\eqref{cond_CK}
contains an $H_n$-factor.
\end{cor}

Write $g=|V(G)|$.
Noting that the condition ``$G$ has no odd components''
implies that $g$ is even.
Considering graphs $G$ with $g$ odd in contrast,
we have the following result.

\begin{thm}\label{thm_OddOrder}
Let $G$ be a connected graph with $g$ odd. Suppose that
\begin{equation}\label{cond_nonempty}
o(G-S)\le 2n|S|,\quad \forall\, S\subseteq V(G),\ S\ne\emptyset.
\end{equation}
Then either $G$ contains an $H_n$-factor, or $G$ is
$H_n^*$-critical.
\end{thm}

\pf Suppose that $G$ contains no $H_n$-factors. Let
\[
(A,B,C,D)
\]
be the $H_n^*$-decomposition of~$G$. From the proof of
Theorem~\ref{thm_main}, we see that $B=\emptyset$, and obtain the
inequality~\eqref{ineq5}.
Together with the condition~\eqref{cond_nonempty},
we see that $A=\emptyset$. Since $G$ is connected, we find that $C=\emptyset$ by
Lemma~\ref{lem_CD}. Hence $G$ is $H_n^*$-critical. This
completes the proof. \qed

We remark that the condition~\eqref{cond_CK} is not necessary
for the existence of an $H_n$-factor in a graph.
Consider the graph
\[
G=K_1+(2n+1)K_{2n+1}
\]
obtained by linking a vertex $K_1$
to all vertices in $2n+1$ copies of the complete graph $K_{2n+1}$. Denote by $C_j$
($1\le j\le 2n+1$) the $j$-th copy of $K_{2n+1}$.
Let $v_j\in V(C_j)$.
Let~$F$ be the factor
consisting of the following $2n+2$ components:
\[
C_1-v_1,\ C_2-v_2,\ \ldots,\ C_{2n}-v_{2n},\ C_{2n+1},\
G\,[\,v_0,\,v_1,\,\ldots,\,v_{2n}\,].
\]
It is easy to verify that~$F$ is an $H_n$-factor.
However,
taking the subset~$S$ to be the single vertex~$v_0$,
we see that the condition~\eqref{cond_CK} does not hold for $G$.

To end this section,
we point out that the coefficient~$2n$
in the condition~\eqref{cond_CK} is a sharp bound
in the sense that for any $\epsilon>0$,
there exists a graph $G$ with a subset $S\subseteq V(G)$ satisfying
\begin{equation}\label{ineq_epsilon}
o(G-S)<(2n+\epsilon)\,|S|,
\end{equation}
and that $G$ contains no $H_n$-factors.
Recall that an $[a,b]$-factor
is a factor $F$ such that $a\le d_F(v)\le b$ for all
vertices $v$.
We need Las Vergnas's theorem~\cite{Las78}.

\begin{thm}[Las Vergnas]\label{thm_Las}
A graph $G$ contains a $[1,n]$-factor
if and only if
for all subsets $S\subseteq V(G)$,
the number of isolated vertices in the subgraph $G-S$
is at most~$n|S|$.
\end{thm}

\begin{thm}\label{thm_epsilon}
For any $\epsilon>0$,
there exists a graph $G$ with a subset $S\subseteq V(G)$ satisfying the $\epsilon$-condition~\eqref{ineq_epsilon}
but with no $H_n$-factors.
\end{thm}

\pf Let $m$ be an integer such that $m>1/\epsilon$.
Let $V_m$ be a set of $m$ isolated vertices,
and $V_{2nm+1}$ a set of $2nm+1$ isolated vertices.
Denote by $K_{m,\,2nm+1}$ the complete bipartite graph
obtained by connecting each vertex in $V_m$
with each vertex in $V_{2nm+1}$.
Setting $S=V_m$ in Theorem~\ref{thm_Las},
we deduce that $K_{m,\,2nm+1}$ contains no $[1,2n]$-factors,
and thus no $H_n$-factors. Moreover,
\[
o(K_{m,\,2nm+1}-V_m)=2nm+1<(2n+\epsilon)\,|V_m|.
\]
This completes the proof. \qed

\section{Sufficient conditions for the existence of an $H_n$-factor}

In this section,
we present some sufficient conditions for
the existence of an $H_n$-factor.
For any vertex~$v$ of~$G$,
denote by $N_G(v)$ the set of neighbors of~$v$.


\begin{thm}\label{thm2}
Let $G$ be a graph without odd components.
If for any non-adjacent vertices $u$ and $v$,
\begin{equation}\label{cond_Neighbor}
|\,N_G(u)\cup N_G(v)\,|>
\max\Bigl\{\,{g-2\over 2n}-1,\ {2g-4\over 4n+1},
\ {g-1\over 2n+1},\ 4n-3\,\Bigr\},
\end{equation}
then $G$ contains an $H_n$-factor.
\end{thm}

\pf Suppose to the contrary that $G$ contains no $H_n$-factors.
By Theorem~\ref{thm_main}, there
exists a nonempty subset~$S\!\subset\!V(G)$ such that
the subgraph $G-S$ has at least $2ns+1$ odd components,
say, $C_1,C_2,\ldots,C_{2ns+1}$, with each $C_i$
has no $H_n$-factors, where $s=|S|$.
Let $c_i=|V(C_i)|$. Suppose that
\begin{equation}\label{order_c}
1\le c_1\le c_2\le\cdots\le c_{2ns+1}.
\end{equation}
It is clear that
\begin{equation}\label{ineq2}
2ns+1\le c_1+c_2+\cdots+c_{2ns+1}\le g-s.
\end{equation}
Therefore
\begin{align*}
c_1&\le{g-s\over 2ns+1},\\[5pt]
c_2&\le{g-s-c_1\over 2ns}.
\end{align*}
It follows that
\begin{align}
c_1+c_2&\le{2(g-s)\over 2ns+1},\label{ineq_c12}\\[5pt]
c_2&\le{g-s-1\over 2ns}.\label{ineq_c2-2}
\end{align}
Moreover, the inequality~\eqref{ineq2} implies that
$s\le s_*$ where
\[
s_*=\frac{g-1}{2n+1}.
\]

Let $u\in V(C_1)$ and $v\in V(C_2)$.
Then
\begin{equation}\label{N}
|\,N_{G}(u)\cup N_{G}(v)\,|
\le s+(c_1-1)+(c_2-1).
\end{equation}
By~\eqref{ineq_c12}, we find that
$|\,N_{G}(u)\cup N_{G}(v)\,|\le h(s)$ where
\[
h(s)={2(g-s)\over 2ns+1}+s-2.
\]
Note that the second derivative $h''(s)>0$.
If $s\ge2$, then we have
\[
|\,N_{G}(u)\cup N_{G}(v)\,|
\le \max\{h(2),h(s_*)\}
=\max\Bigl\{\,\frac{2g-4}{4n+1},\ \frac{g-1}{2n+1}\,\Bigr\},
\]
contradicting to the condition~\eqref{cond_Neighbor}.
Otherwise $s=1$.
In this case,
if $c_2\le 2n-1$,
then $c_1\le 2n-1$ by~\eqref{order_c}.
By~\eqref{N}, we have
\[
|\,N_{G}(u)\cup N_{G}(v)\,|
\le c_1+c_2-1
\le 4n-3,
\]
contradicting to~\eqref{cond_Neighbor}.
So $c_2\ge 2n$.
It is easy to verify that any complete graph~$K_m$
with $m\ge 2n$ has an $H_n$-factor.
Since~$C_2$ contains no $H_n$-factors,
we deduce that~$C_2$ is not complete.
So there exist vertices~$u'$ and~$v'$ which are not
adjacent in~$C_2$. By~\eqref{ineq_c2-2}, we have
\[
|\,N_{G}(u')\cup N_{G}(v')\,|\le s+c_2-2\le\frac{g-2}{2n}-1,
\]
contradicting to~\eqref{cond_Neighbor}.
This completes the proof. \qed

Observe that when $g\ge 8n^2+2n+2$, one has
\[
\max\Bigl\{{g-2\over 2n}-1,\,{2g-4\over 4n+1},\,{g-1\over 2n+1},\,4n-3\Bigr\}
={g-2\over 2n}-1.
\]
This results in the following corollary immediately.

\begin{cor}
Let $G$ be a graph without odd components. If
$g\ge 8n^2+2n+2$,
and for any non-adjacent vertices $u$ and $v$,
\[
|\,N_G(u)\cup N_G(v)\,|>{g-2\over 2n}-1,
\]
then $G$ contains an $H_n$-factor.
\end{cor}

Now we give another sufficient condition
for the existence of an $H_n$-factor of a graph.
A graph~$G$ is said to be $k$-connected
if it is connected when
fewer than~$k$~vertices are removed from~$G$.
Let~$u$ and~$v$ be non-adjacent vertices of $G$.
Denote by $G+uv$ the graph obtained by adding
the edge $(u,v)$ to $G$.

%

\begin{thm}
Let $G$ be a $k$-connected simple graph with~$g$ even.
Let~$u$ and~$v$ be non-adjacent vertices of~$G$
such that
\begin{equation}\label{cond-kconnected}
|\,N_G(u)\cup N_G(v)\,|\ge g-2nk.
\end{equation}
Then $G$ has an $H_n$-factor if and only if
the graph $G+uv$ has an $H_n$-factor.
Moreover, the lower bound $g-2nk$ in~\eqref{cond-kconnected}
is best possible.
\end{thm}

\pf The necessity is obvious. We shall prove the sufficiency.
Suppose to the contrary that $G$ has no $H_n$-factors.
Since $G$ is a connected graph with $g$ even,
we deduce that $G$ has no odd components.
By Theorem~\ref{thm_main}, there exists a nonempty subset~$S\subset V(G)$ such that
\begin{equation}\label{eq2}
o(G-S)\ge2ns+1,
\end{equation}
where $s=|S|$.
Since $G$ is $k$-connected, it is easy to see that $s\ge k$.
Let $C_1,C_2,\ldots,C_q$ be the components of $G-S$, and $c_i=|V(C_i)|$.
Then $q\ge 2n+1$ by~\eqref{eq2}.
By Theorem~\ref{thm_main},
we can suppose that $c_i$ is odd and $C_i$ has no $H_n$-factors
for any $1\le i\le 2ns+1$.

Let $u$ and $v$ be non-adjacent vertices. We have three cases.
\begin{itemize}
\item[(i)]
$u$ and $v$ belong to the same component of $G-S$.
In this case, we can suppose that $u,v\in C_i$ for some $1\le i\le q$.
By~\eqref{eq2}, we have
\[
|\,N_G(u)\cup N_G(v)\,|
\le s+(c_i-2)
=g-\sum_{j\ne i}c_j-2
\le g-2ns-2
\le g-2nk-2.
\]
\item[(ii)]
$u$ and $v$ belong to distinct components of $G-S$.
In this case, we can suppose that $u\in C_i$ and $v\in C_j$
where $1\le i<j\le q$.
By~\eqref{eq2}, we have
\begin{align*}
|N_G(u)\cup N_G(v)|\le s+c_i+c_j-2
=g-\sum_{h\not\in\{i,\,j\}}c_h-2
\le g-2nk-1.
\end{align*}
\item[(iii)]
One of $u$ and $v$ belongs to the set $S$.
Let $F$ be an $H_n$-factor of $G+uv$.
Let $m$ be the total degree of vertices of~$S$ in~$F$, i.e.,
\[
m=\sum_{v\in S}d_F(v).
\]
Since the component $C_i$ contains no $H_n$-factors,
and the edge $(u,v)$ is not contained in $G-S$,
there exists an edge of~$F$ connecting one
vertex in $C_i$ and another vertex in~$S$.
Therefore, each~$C_i$ corresponds
an edge with one end in~$S$.
It follows that $m\ge 2ns+1$.
On the other hand, each vertex in~$S$ has degree at
most~$2n$ in~$F$. So $m\le 2ns$, a contradiction.
This proves the sufficiency.
\end{itemize}

Now we shall prove that the bound $g-2nk$
is best possible. Let $k$ be an odd number, and
\[
G=K_{k}+(2nk+1)K_1.
\]
Then $g=2nk+k+1$ is even.
It is easy to check that $G$ is $k$-connected and
\[
|N_G(u)\cup N_G(v)|=g-2nk-1
\]
for any non-adjacent vertices $u,v\in V(G)$.
It suffices to show that $G+uv$ contains an $H_n$-factor
while $G$ does not.

Denote by $u_1,u_2,\ldots,u_k$ the vertices of the subgraph~$K_k$,
and by $v_1,v_2,\ldots,v_{2nk+1}$ the remaining vertices in~$G$.
Suppose to the contrary that~$G$ has an $H_n$-factor~$F$.
Then the degree of each vertex $v_i$ is at least 1 in~$F$.
Note that the neighbor of $v_i$ must be some $u_j$.
So there exists some $u_j$ of degree at least
$2n+1$ in~$F$. It follows that the degree of $u_j$ does not
belong to the set $H_n$, a contradiction.
Hence~$G$ has no $H_n$-factors.

Now we shall show that $G+uv$ contains an $H_n$-factor
for any non-adjacent vertices $u$ and $v$.
In fact, since each vertex $u_i$ is saturated,
we can suppose without loss of generality that
the non-adjacent vertex pair $(u,v)$
is taken to be $(v_{2nk},v_{2nk+1})$.
A factor $F$ of $G+uv$ consists of the edge~$uv$
and the following $k$ components:
\begin{align*}
&G[u_1,\,v_1,\,v_2,\,\ldots,\,v_{2n}],\\[5pt]
&G[u_2,\,v_{2n+1},\,v_{2n+2},\,\ldots,\,v_{4n}],\\[5pt]
&\quad\vdots\\[5pt]
&G[u_k,\,v_{2n(k-1)+1},\,v_{2n(k-1)+2},\,\ldots,\,v_{2nk-1}].
\end{align*}
It is straightforward to verify that $F$ is an $H_n$-factor.
This completes the proof. \qed

Let $H$ be an allowed prescription.
In~\cite[Lemma~(3.5)]{Lov72},
Lov\'asz gave the following result
describing the $H$-decomposition of a graph
when a vertex in $A_H$ is removed.

\begin{lem}[Lov\'asz]\label{lem_v}
Let $(A,B,C,D)$ be the $H$-decomposition of~$G$.
Let $v$ be a vertex in~$A$, and
$(A',B',C',D')$ the $H$-decomposition of the subgraph $G-v$.
Then
\[
A'=A-v,\quad
B'=B,\quad
C'=C,\quad
D'=D.
\]
\end{lem}

\begin{thm}
Let $G$ be a graph without odd components.
Then $G$ contains an $H_n$-factor
if the subgraph $G-v$ contains an $H_n$-factor
for all vertices $v$.
\end{thm}

\pf Suppose that $G$ has no $H_n$-factors.
Let $(A,B,C,D)$
be the $H_n^*$-decomposition of~$G$.
From the proof of Theorem~\ref{thm_main},
we see that $B=\emptyset$
and
\begin{equation}\label{eq3}
2n|A|<c(D).
\end{equation}
Moreover, every component of $G[D]$ is odd.

Assume that $A\ne\emptyset$.
Let $v\in A$,
and $(A',B',C',D')$ be the $H_n^*$-decomposition of $G-v$.
By Lemma~\ref{lem_Lovasz1}, we have
\begin{equation}\label{ineq4}
c(D')\le 2n|A'|.
\end{equation}
By~\eqref{eq3}, \eqref{ineq4} and Lemma~\ref{lem_v},
we deduce that
\[
2n|A|\le c(D)-1=c(D')-1\le 2n|A'|-1=2n(\,|A|-1\,)-1,
\]
a contradiction. So $A=\emptyset$.
By Lemma~\ref{lem_CD},
any component of $G[D]$
is an odd component of~$G$.
But $G$ has no odd components, a contradiction.
This completes the proof.
\qed

\end{document}